\documentclass{article}
\usepackage[utf8]{inputenc}
\usepackage{amsmath}
\usepackage{amssymb}
\title{Review on rationality problems of algebraic k-tori}
\author{Youngjin Bae }
\newtheorem{defn}{Definition}[section]
\newtheorem{ex}{Example}[section]
\newtheorem{thm}{Theorem}[section]
\date{}

\begin{document}

\maketitle

\begin{abstract}
   Rationality problems of algebraic $k-tori$ are closely related to rationality problems of the invariant field, also known as Noether's Problem. 
   We describe how a function field of algebraic $k-tori$ can be identified as an invariant field under a group action and that a $k-tori$ is rational if and only if its function field is rational over $k$. We also introduce character group of $k-tori$ and numerical approach to determine rationality of $k-tori$.  
\end{abstract}
\vspace{7mm}
\tableofcontents
\newpage
\section{Introduction}

\vspace{5mm} 

Let $k$ be a field and $K$ is a finitely generated field extension of $k$. $K$ is called \textit{rational over k} or \textit{k-rational} if $K$ is isomorphic to $k(x_{1},...,x_{n})$ where $x_{i}$ are transcendental over $k$ and algebraically independent. There are also relaxed notions of rationality. $K$ is called \textit{stably k-rational} if $K(y_1,...,y_m)$ is $k-rational$ for some transcendental and algebraically independent $y_{i}$. $K$ is called $k-unirational$ if $k \subset K \subset k(x_1,...,x_n)$ for some pure transcendental extension $k(x_1,...,x_n)/k$.

The Noether's Problem is the question of rationality of the invariant field under finite group action. For example, if $K=\mathbb{Q}(x_1,x_2)$ and $G=\{1,\sigma\}\cong C_2$ and $G$ acts on $K$ as permutation of variables $x_1,x_2$ (i.e. $\sigma$ fixes $\mathbb{Q}$, $\sigma(x_1)=x_2$ and $\sigma(x_2)=x_1$), then the invariant field $K^G$ is $\mathbb{Q}-rational$.
\begin{ex}
$K=\mathbb{Q}(x,y)$ and $G \cong C_2$, acting on $K$ as permutation of variables. Let $\frac{f}{g} \in K^G$, $f,g$ are coprime. We have $$\frac{f(x,y)}{g(x,y)}=\sigma(\frac{f(x,y)}{g(x,y)})=\frac{f(y,x)}{g(y,x)}$$ 
By observing that $gcd(f(x,y),g(x,y))=gcd(f(y,x),g(y,x))=1$, we have $f(x,y)=f(y,x)$ and $g(x,y)=g(y,x)$. \newline Therefore, $K^G=\{\frac{f(x,y)}{g(x,y)}|f,g$ are symmetric$\}$, field of fractions (quotient field) of $S=\{f \in \mathbb{Q}[x,y]|f(x,y)=f(y,x)\}$. 
It is easy to see that $\psi:S \to \mathbb{Q}[s,t]$ is isomorphism, where 
$$\psi(x+y)=s, \hspace{3mm} \psi(xy)=t$$
Therefore, $S \cong \mathbb{Q}[x,y]$ and $K^G \cong \mathbb{Q}(x,y)$, $\mathbb{Q}-rational$.
\end{ex}
\vspace{3mm}
We can also consider case of $G$ acting on both of coefficients and variables. 
\vspace{3mm}
\begin{ex}
$K=\mathbb{C}(x,y)$ and $G=Gal(\mathbb{C}/\mathbb{R})=\{1,\sigma\}\cong C_2$. Suppose $G$ acts on $K$ by permuting $x,y$ and as complex conjugation on coefficients.\newline
For example, $\sigma(ix+(1-i)xy+y^2)=-iy+(1+i)yx+x^2$.
Then, $K^G \cong \mathbb{R}(x,y)$, is $\mathbb{R}-rational$.

\textbf{Proof.} For $\frac{f(z,w)}{g(z,w)} \in K^G$, where $f,g$ are coprime, $\sigma(f)$ and $\sigma(g)$ are also coprime. From $\frac{f}{g}=\frac{\sigma(f)}{\sigma(g)}$, we have $f=\sigma(f)$ and $g=\sigma(g)$. Thus, $K^G$ is quotient field of $S$ where $S:=\{f(z,w) \in \mathbb{C}[z,w]|f=\sigma(f)\}$.

Define a map $\psi : S \to \mathbb{R}[x,y]$ as 
$$z=x+yi, w=x-yi$$
$$and$$
$$\psi(f)(x,y)=f(z,w)$$

The coefficients of $\psi(f)$ are real numbers. This is because, if we let $f(z,w)=\sum_{n,m}a_{n,m}z^{n}w^{m}$, we have that $$\psi(f)(x,y)=f(z,w)=\sigma(f(z,w))=\sigma(\sum_{n,m}a_{n,m}z^{n}w^{m})=\sum_{n,m}\overline{a_{n,m}}w^{n}z^{m}$$
$=\sum_{n,m}\overline{a_{n,m}(x+iy)^{n}(x-iy)^{m}}=\overline{\psi(f)(x,y)}$.

\hspace{3mm}

Therefore, $\psi(f)=\overline{\psi(f)}$, $\psi(f) \in \mathbb{R}[x,y]$. It is easy to see that $\psi$ is actually isomorphism, $S \cong \mathbb{R}[x,y]$, and $K^G \cong \mathbb{R}(x,y)$. 

\label{exa}
\end{ex}

Another perspective to view this \textit{change of variables} is identifying the field with rational function field of algebraic $k-tori$. (see \textbf{Example \ref{exc}} and \textbf{Example \ref{exd}})

\section{Algebraic $k-tori$}

Let $k$ be a field. Then $\mathbb{A}_{k}^{n}$ is \textit{n-dimension affine space}  over the field $k$, simply $k^n$ with usual vector space structure on it.  A subset $X$ of $\mathbb{A}_{k}^{n}$ is an \textit{algebraic k-variety} (\textit{k-variety} in short) if it is a set of zeros of a system of equations with $n$ variables $x_{1},...x_{n}$ over $k$. The ideal of polynomials that vanish on every points of $X$ will be denoted by $I(X)$. The \textit{coordinate ring} of a variety $X$ is defined to be the quotient 
$$A(X):=k[x_1,...,x_n]/I(X)$$
Projective varieties can be similarly defined as the set of zeros of a system of homogeneous equations. \textit{Projective $n-space$} $\mathbb{P}_{k}^n$ is defined as set of lines passing the origin in $\mathbb{A}_{k}^{n+1}$.

If $X,Y$ are varieties, a map $f: X \to Y$ is called \textit{regular} if it can be presented as fraction of polynomials $p/q$, where $q$ does not vanishes in $X$. A map $f: X \to Y $ is called \textit{rational} if it is regular on Zariski open dense set. (Formally, a regular map is defined as an equivalence class of pairs $<U,f_U>$ where $U$ is Zariski open subset of $U$. See  \cite{H77})
Let $X$ be a variety, $K(X)$ is the \textit{rational function field}, or \textit{function field} in short, the set of rational maps $f:X \to \mathbb{A}_{k}$. For example, if $X$ is an affine variety over algebraically closed field $k$, $K(X)$ is quotient field of $A(X)$. 

\begin{ex}
Let $X=\{(x,y) \in \mathbb{A}_{\mathbb{C}}^2|xy=1\}$ be a variety over $\mathbb{C}$. \newline Then, $A(X)=\mathbb{C}[x,y]/(xy-1)\cong\mathbb{C}[x,\frac{1}{x}]$ and $K(X)\cong\mathbb{C}(x)$.
\end{ex}

Two varieties $X,Y$ are \textit{isomorphic} (resp. \textit{birationally isomorphic}) if there is a bijective regular map (resp. rational map) $f:X \to Y$ and its inverse is also regular (resp. rational). 

A variety $X$ in $\mathbb{A}_{k}^{n}$ is an \textit{algebraic group} if it has a group structure on it, where the group operation and inversions are regular maps. (i.e. $\ast : X \times X \to X$ and $^{-1} : X \to X$ are regular) 

Algebraic $k-tori$, or algebraic $k-torus$, is a special type of algebraic group over $k$. We call an algebraic group as $k-tori$ when it is isomorphic to some power of multiplicative group over $\overline{k}$, the algebraic closure of $k$. 

\begin{defn}[Multiplicative Group]
Let $k$ be a field, the multiplicative group $\mathbb{G}_{m}(k)$ is algebraic group in $\mathbb{A}_{k}^{2}$, defined as $\{(x,y)\in\mathbb{A}_{k}^{2}|xy=1\}$, with operation $\cdot: \mathbb{G}_{m}(k)\times \mathbb{G}_{m}(k) \to \mathbb{G}_{m}(k)$ of $(x,\frac{1}{x})\cdot(y,\frac{1}{y})=(xy,\frac{1}{xy})$
\end{defn}

\begin{ex}
$\mathbb{G}_{m}(\mathbb{R})$ is the curve $xy=1$ on the real affine plane. It is isomorphic to $\mathbb{R}^{\times}$ as a group. $((x,y) \to x$ is group isomorphism.$)$
\end{ex}

As field changes, same system of equations can define different varieties. For instance, the equation $xy=1$ in previous example defines $\mathbb{G}_{m}(\mathbb{C})$ in $\mathbb{A}_{\mathbb{C}}^{2}$, which is different from $\mathbb{G}_{m}(\mathbb{R})$. If $E$ is a field and $F$ is its algebraic closure, an irreducible variety $V$ over $F$ entails the ring of equations, $I$. If $I$ happens to be in $E[\textbf{x}]$ (ring of polynomials over $E$), we can define $V(E)$, a variety over $E$ defined by equations in $I$. This can be viewed as \textit{restriction} of scalar. Extension of scalar can be defined similarly.

\begin{defn}[Algebraic $k$-tori] 
Let $k$ be a field with algebraic closure $\overline{k}$. If $T$ is an algebraic group over $k$, it is $k-torus$ if and only if 

\begin{center}$T(\overline{k}) \cong (\mathbb{G}_{m}(\overline{k}))^{r}$ 
\end{center}
for some $r$. The $r$ is called dimension of $T$.
\end{defn}

\begin{ex}
$T=\mathbb{G}_{m}(\mathbb{R})$ is one dimensional $\mathbb{R}-tori$. This is because $T(\mathbb{C})=\mathbb{G}_{m}(\mathbb{C})$.
\end{ex}

From now, let $k^{\times}=\mathbb{G}_{m}(k)$ be the one dimensional torus over $k$.
There are two one-dimensional  $\mathbb{R}$-tori, one can be recognized as $\mathbb{R}^{\times}$, the other one can be recognized as $SO(2)$ as a group. 

\begin{ex}
The norm one torus $N$ is a real algebraic group in $\mathbb{A}_{\mathbb{R}}^{2}$, defined by equation $x_{1}^2+x_{2}^2=1$ $(i.e.$   $ N=\{(x_1,x_2) \in \mathbb{A}_{\mathbb{R}}^{2}|x_{1}^2+x_{2}^2=1\})$, and operation $\cdot:N \times N \to N$ such that $$(x_{1},x_{2})\cdot(y_{1},y_{2})=(x_{1}y_{1}-x_{2}y_{2},
x_{1}y_{2}+x_{2}y_{1})$$
Indeed, $N$ is isomorphic to $SO(2)$ as a group. \newline
Also, $N(\mathbb{C})=\{(x_1,x_2) \in \mathbb{A}_{\mathbb{C}}^{2}|x_{1}^2+x_{2}^2=1\}$ is isomorphic to $C^{\times}$ as algebraic group. The map $\psi : N(\mathbb{C}) \to \mathbb{C}^{\times} $ 
$$\psi(x_1,x_2)=x_1+ix_2 $$
is isomorphism. 
Therefore, $N$ is one dimensional real torus.
\end{ex}

\vspace{3mm}

If $T$ is a $k-torus$, $T$ is called \textit{split over K} if it satisfies $T(K) \cong (K^{\times})^{s}$ for some extension $K/k$ and some $s$. For instance, $\mathbb{R}^{\times}$ is split over $\mathbb{R}$, $N$ is not. \newline
It is easy to find split torus such as $(\mathbb{R}^{\times})^{2}$ or $(\mathbb{R}^{\times})^{3}$, being another torus. Also, for any integer $r$, $N^{r}$ is $r$-dimensional $\mathbb{R}-tori$. Meanwhile, there are also some non-trivial(not a product of low-dimensional torus) torus.  

\begin{ex}
Let $P$ be a real algebraic group in $\mathbb{A}_{\mathbb{R}}^{4}$, defined as $$P=\{(x_{1},x_{2},x_{3},x_{4}) \in \mathbb{A}_{\mathbb{R}}^{4}|x_{1}x_{3}-x_{2}x_{4}=1  , x_{1}x_{4}+x_{2}x_{3}=0\}$$
Alternatively, $$P=\{A \in M_{2\times2}(\mathbb{R}) \hspace{1mm} | \hspace{1mm} AA^{t}=
\quad
\begin{pmatrix} 
s & 0 \\
0 & s^{-1} 
\end{pmatrix}
\quad
s \in \mathbb{R}\backslash\{0\}\}$$
and operation $\cdot:P \times P \to P$ such that 
$$(x_{1},x_{2},x_{3},x_{4})\cdot(y_{1},y_{2},y_{3},y_{4})=(x_{1}y_{1}-x_{2}y_{2},x_{1}y_{2}+x_{2}y_{1},x_{3}y_{3}-x_{4}y_{4},x_{3}y_{4}+x_{4}y_{3})$$
Which is compatible with complex multiplication of 
$$(x_{1}+x_{2}i,x_{3}+x_{4}i)\cdot(y_{1}+y_{2}i,y_{3}+y_{4}i)$$
Moreover, $P(\mathbb{C})$ is isomorphic to $(\mathbb{C}^{\times})^{2}$, by sending $$(x_{1},x_{2},x_{3},x_{4}) \to ((x_{1}+x_{2}i,x_{3}+x_{4}i),(x_{1}-x_{2}i,x_{3}-x_{4}i))=((z,\frac{1}{z}),(w,\frac{1}{w}))$$Therefore, $P$ is 2-dimensional $\mathbb{R}-tori$.

\label{exc}
\end{ex} 
\vspace{3mm}
By tracking the function fields of $P(\mathbb{R})$ and $P(\mathbb{C})$, we have the same trick of change of variables as in \textbf{Example \ref{exa}}.

\vspace{3mm}

\begin{ex}
In the previous example, the coordinate ring of $P(\mathbb{C})$ is
$$A(P(\mathbb{C}))=\mathbb{C}[x_1,x_2,x_3,x_4]/(x_{1}x_{3}-x_{2}x_{4}-1,x_{1}x_{4}+x_{2}x_{3}) \cong \mathbb{C}[z,\frac{1}{z},w,\frac{1}{w}]$$
where $z=x_1+x_2i$ and $w=x_1-x_2i$. The function field of $P(\mathbb{C})$ is 
$$K(P(\mathbb{C}))\cong \mathbb{C}(z,w)$$
Let $G=Gal(\mathbb{C}/\mathbb{R})$ acts on $K(P(\mathbb{C}))$ as in \textbf{Example \ref{exa}}. Observe that the coordinate ring of $P(\mathbb{R})$ is $A(P(\mathbb{R}))=A(P(\mathbb{C}))^G$ and the function field of $P(\mathbb{R})$ is $K(P(\mathbb{R}))=K(P(\mathbb{C}))^G \cong \mathbb{C}(z,w)^G$ (note that $G$ actions on $K(P(\mathbb{C}))$ and $\mathbb{C}(z,w)$ are equivalent through the isomorphism). In short, we have that 
$$K(P(\mathbb{R}))\cong \mathbb{C}(z,w)^G$$
Therefore, when $G=Gal(\mathbb{C}/\mathbb{R})$ action on $C(z,w)$ is given, we can convert the rationality problem to the rationality problem of $K(P(\mathbb{R}))$, the function field of $P(\mathbb{R})$. In this sense, the following definition and theorem are natural. 
\label{exd}
\end{ex}

\begin{defn}[Rationality of $k-variety$] 
We say that a variety $X$ over $k$ is rational if, equivalently, \newline

(1) $X$ is birationally isomorphic to $\mathbb{P}_{k}^{n}$ for some $n$.

(2) $K(X)\cong k(x_1,..,x_n)$
\end{defn}

\vspace{3mm} 
If $K/k$ is Galois extension,  a $k-tori$ $T$ is $K-rational$ if it is rational as a $K$-variety $T(K)$. If $k$ is algebraically closed, there is unique $n$-dimension tori $T_n=(k^{\times})^{n}$. Since the function field of $T_n$ is $k(x_1,...,x_n)$, thus $T_n$ is $k$-rational.

\vspace{3mm}

\begin{thm}
The following two problems are equivalent. 
\newline

(1) The rationality problem of $n$ dimensional $k-tori$ $T$

(2) The rationality problem of invariant field $K^G$ \newline

\noindent where $G=Gal(\overline{k}/k)$ and $K=k(x_1,...,x_n)$.
\end{thm}

\vspace{6mm}

There is a connection between the $G$ action on $K$ and $k-tori$ $T$, connecting the two rationality problems given in the previous theorem. To be specific, the character group of $T$ determines both the $G$ action and $T$ uniquely.

\section{Character group of $k-tori$}

\begin{defn}[Character group of $k-tori$]
Let $T$ be $k-tori$. Then $\mathbb{X}(T)$, the \textit{character group} of $T$ is the set of algebraic group homomorphisms(a regular map preserving the group structure) from $T$ to $\overline{k}^{\times}$, denoted by $Hom(T,\mathbb{G}_{m})$ or $Hom(T,\overline{k}^{\times})$. 
\end{defn}

The character group $\mathbb{X}(T)$ of $T$ has a group structure defined by component-wise multiplication. Also, if $T$ is split over $L$ for finite Galois extension of base field $k$, $G=Gal(L/k)$ acts on $\mathbb{X}(T)$. Moreover, it is known that $\mathbb{X}(T)$ is torsion-free $\mathbb{Z}$-module(i.e. isomorphic to $\mathbb{Z}^{n}$ for some $n$). Therefore, $\mathbb{X}(T)$ is a $G-lattice$ (a free $\mathbb{Z}-module$ with $G$-action).

\begin{ex}
If $T=\mathbb{C}^{\times}$ is multiplicative group of $\mathbb{C}$, then $\mathbb{X}(T)$ is set of regular functions $f: \mathbb{C}^{\times} \to \mathbb{C}^{\times}$ such that $f(xy)=f(x)f(y)$ for $x,y \in \mathbb{C}^{\times}$. Since $f$ is a rational function, it is a meromorphic function over $\mathbb{C}$. Also, we have $f(\mathbb{C}^{\times}) \subset \mathbb{C}^{\times}$, which implies 0 is the only point where $f$ can have zeros or poles. Therefore, $f(t)=t^{n}$ for some $n \in \mathbb{Z}$. If we write a function $t \to t^{n}$ as $t^n$, we have $$\mathbb{X}(T)=\{t^n|n \in \mathbb{Z}\} \cong \mathbb{Z}^{1}$$ as a group. 
$G=Gal(\mathbb{C}/\mathbb{C})=\{id\}$ acts trivially on $\mathbb{X}(T)$.
\end{ex}

\vspace{3mm}

\noindent In general, if $k$ is algebraically closed, the character group of $(k^{\times})^{n}=\mathbb{G}_{m}^{n}$ is

\vspace{2mm}

\noindent $\mathbb{X}(\mathbb{G}_{m}^{n})=\{f_{t_1,...t_n}:\mathbb{G}_{m}^{n} \to \mathbb{G}_{m}| f_{t_1,...t_n}(x_1,...x_n)=\prod_{i}x_{i}^{t_i}, t_i \in \mathbb{Z} \}$

\vspace{2mm}

\noindent $=\prod_{i=1}^{n} \{f_{t}: \mathbb{G}_{m} \to \mathbb{G}_{m} | f_{t}(x_i)=x_{i}^{t}, t \in \mathbb{Z}\} \cong \mathbb{Z}^n$

\vspace{3mm}

\begin{ex}
Let $P$ be the 2-dimension $\mathbb{R}-tori$ in \textbf{Example \ref{exc}}. Then, the character group of $P$ is 
$$\mathbb{X}(P)=\{f_{t_1,t_2}:P \to \mathbb{C}^{\times}|f_{t_1,t_2}(x_1,x_2,x_3,x_4)=(x_1+x_{2}i)^{t_1}(x_1-x_{2}i)^{t_2}\}$$

\noindent Let $z=x_1+x_{2}i$, $w=x_1-x_{2}i$, then we have the natural extension of $\mathbb{X}(P)$ to $\mathbb{X}(P(\mathbb{C}))$ 

$$\mathbb{X}(P(\mathbb{C}))=\{f_{t_1,t_2}:P(\mathbb{C}) \to \mathbb{C}^{\times}|f_{t_1,t_2}((z,\frac{1}{z}),(w,\frac{1}{w}))=z^{t_1}w^{t_2}\}\cong \mathbb{Z}^2$$

\noindent Observe that the complex conjugation $\sigma \in G$, exchanges $z$ and $w$, thus acting on $\mathbb{Z}^2$ as $2 \times 2$ matrix $
\begin{bmatrix} 
0 & 1 \\
1 & 0 
\end{bmatrix}
$.

\label{exe}
\end{ex}

\vspace{7mm}

\noindent It is known that when a $G=Gal(K/k)$ action (as $\mathbb{Z}$-linear function) on $\mathbb{Z}^{n}$ is given, there exists unique  $n$-dimensional $k-tori$ which has the given $G-lattice$ as its character group. Furthermore, there are conditions of $G-lattice$ corresponding to the rationality conditions of $k-tori$ and of invariant fields. 

\vspace{5mm}

\section{Flabby resolution and numerical approach}

This section contains many results in \cite{AA17}. Let $G$ be a group and $M$ be a $G-lattice$ ($M \cong \mathbb{Z}^{n}$ as group and has $G$-linear action on it). $M$ is called a \textit{permutation G-lattice} if $M\cong \bigoplus_{1\leq i \leq m}\mathbb{Z}[G/H_i]$ for some subgroups $H_1,...,H_m$ of $G$ (equivalently, there exists a $\mathbb{Z}$-basis of $M$ such that $G$ acts on $M$ as permutation of the basis). $M$ is called \textit{stably permutation G-lattice} if $M\bigoplus P \cong Q$ for some permutation $G-lattices$ $P$ and $Q$. $M$ is called \textit{invertible} if it is a direct summand of a permutation $G$-lattice, i.e. $P\cong M \bigoplus M'$ for some permutation $G$-lattice $P$ and $M'$. 

\begin{defn}[1st Group Cohomology]
Let $G$ be a group and $M$ be a $G$-lattice. For $g \in G$ and $m \in M$, let $g.m=m^g$ be $g$ acting on $m$. The first group cohomology $H^{1}(G,M)$ is a group defined as 

$$H^{1}(G,M)=Z^{1}(G,M)/B^{1}(G,M)$$

\vspace{3mm}

\noindent where $Z^{1}(G,M)=\{f:G \to M | f(gh)=f(g)^{h}f(h)\}$ and $B^{1}(G,M)=\{f:G \to M|f(g)=m_{f}^{g}m_{f}^{-1}\; for \; some \; m_{f} \in M\}$
\end{defn}

\newpage

\noindent $H^{1}(G,M)=0$ simply implies that if $f: G \to M$ satisfies $f(gh)=f(g)^{h}f(h)$, then there exists $m \in M$ such that $f(g)=m^{g}m^{-1}$. $M$ is called \textit{coflabby} if $H^{1}(G,M)=0$. 

\begin{defn}[-1st Tate Cohomology]
Let $G$ be finite group of order n and $M$ be a $G$-lattice. The -1st group cohomology $\hat{H}^{-1}(G,M)$ is a group defined as 

$$\hat{H}^{-1}(G,M)=Z^{-1}(G,M)/B^{-1}(G,M)$$

\vspace{3mm}

\noindent where

$$Z^{-1}(G,M)=\{m \in M| \sum_{g \in G}m^g=0\}$$,
$$B^{-1}(G,M)=\{\sum_{g \in G}m_{g}^{g-id}|m_{g} \in M\}$$
\end{defn}

\noindent Similarly, $M$ is called \textit{flabby} if $\hat{H}^{-1}(G,M)=0$. It is clear that a $k-tori$ is rational if and only if $\mathbb{X}(T)$ is permutation $G$-lattice. Thus, the rationality problems of $k-tori$ and invariant fields can be reduced into problem of finding permutation $G$-lattice(equivalent to find finite subgroup of $GL(n,\mathbb{Z})$. However, this problem is not solved yet, even though there are many results in weakened problems. 

Let $C(G)$ be the category of all $G$-lattices and $S(G)$ be the category of all permutation $G$-lattices. Define equivalence relation $~$ on $C(G)$ by $M_1~M_2$ if and only if there exist $P_1,P_2 \in S(G)$ such that $M_1\bigoplus P_1\cong M_2\bigoplus P_2$. Let $[M]$ be equivalence class containing $M$ under this relation. 

\begin{thm}\textnormal{(Endo and Miyata \cite[Lemma 1.1]{EM75}, Colliot-Th\'el\`ene and Sansuc \cite[Lemma 3]{CTS77})} 
For any $G$-lattice $M$, there is a short exact sequence of $G$-lattices $0 \to M \to P \to F \to 0$ where $P$ is permutation and $F$ is flabby. 
\end{thm}

In the previous theorem, $[F]$ is called the \textit{flabby class} of $M$, denoted by $[M]^{fl}$. 

\begin{thm}\textnormal{(Akinari and Aiichi \cite[17pp]{AA17})} If $M$ is stably permutation, then $[M]^{fl}$. If $M$ is invertible, $[M]^{fl}$ is invertible. 
\end{thm}

\noindent It is not difficult to see that 
$$M \; is \; permutation \; \Rightarrow \; M \; is \; stably \; permutation$$

\noindent Furthermore, it is true that 
$$M \; is \; stably \; permutation \; \Rightarrow \; M \; is \; invertible \; \Rightarrow M \; is \; flabby \; and \; coflabby  $$

\noindent In \cite{AA17}, they gave the complete list of stably permutation lattices for dimension 4 and 5 by computing $[M]^{fl}$ for finite subgroup of $GL(n,\mathbb{Z})$, which is equivalent to classifying stably rational tori. Thus, the rationality problems for low dimensional $k-tori$ can be resolved by finding conditions which can determine a stably permutation $M$ is permutation or not. 

\newpage


\end{document}